\documentclass{amsart}
\usepackage{amssymb,amsmath,amsthm}

\theoremstyle{plain}
\newtheorem{lemma}{Lemma}[section]
\newtheorem{theorem}[lemma]{Theorem}
\newtheorem{corollary}[lemma]{Corollary}
\newtheorem{sublemma}[lemma]{Sublemma}

\theoremstyle{definition}
\newtheorem{example}[lemma]{Example}

\newcommand{\K}{\mathbf{k}}
\newcommand{\f}{\varphi}
\newcommand{\tp}{\overline}

\DeclareMathOperator{\gr}{gr}
\DeclareMathOperator{\ak}{AK}
\DeclareMathOperator{\trdeg}{trdeg}
\DeclareMathOperator{\grdeg}{grdeg}
\DeclareMathOperator{\Spec}{Spec}
\DeclareMathOperator{\EXP}{EXP}
\DeclareMathOperator{\ch}{char}
\DeclareMathOperator{\Frac}{Frac}
\begin{document}
\title{Return of $x + x^2y + z^2 + t^3 = 0$}
\author{A. Crachiola}
\address{Department of Mathematics\\
Wayne State University\\
Detroit, MI 48202\\
USA\\} \email{ crach@math.wayne.edu}

\date{February 18, 2004}

\begin{abstract}
We develop techniques for computing the AK invariant of domains
with arbitrary characteristic. As an example, we show that for any
field $\K$ the ring $\K[X,Y,Z,T] / (X + X^2 Y + Z^2 + T^3)$ is not
isomorphic to a polynomial ring over $\K$.
\end{abstract}
\maketitle

\section{Introduction}

All rings in this paper are commutative with identity. Throughout
this paper, let $\K$ denote a field of arbitrary characteristic,
and let $\K^* = \K \setminus 0$. For a ring $A$, let $A^{[n]}$
denote the polynomial ring in $n$ indeterminates over $A$. Let
$\mathbf{C}$ denote the complex number field.

The AK invariant was introduced by L. Makar-Limanov \cite{ML1} to
show that P. Russell's \cite{KorRus} threefold
\[
x + x^2 y + z^2 + t^3 = 0
\]
over $\mathbf{C}$ is not isomorphic to $\mathbf{C}^3$. The
invariant was defined in that context as the intersection of the
kernels of all locally nilpotent derivations on the coordinate
ring of a variety. Since this first application, the AK invariant
has been applied by several people (for instance
\cite{bandman,dubouloz,KKMLR,ML2,villa}), mainly in the realm of
affine algebraic geometry. To the author's knowledge, this work
has all been conducted under the restriction of zero
characteristic. Superficially, this is reasonable, since
derivations don't behave as nicely on rings with prime
characteristic $p$. The kernel of such a derivation doesn't convey
the appropriate information in this setting, since it will always
contain the $p$th power of every element.

Nevertheless, the apparent advantage to restricting the
characteristic may be only a matter of perception. In the zero
characteristic situation, locally nilpotent derivations on the
coordinate ring correspond with algebraic additive group actions
on the variety, and, unlike locally nilpotent derivations, these
actions maintain their attractive properties in the prime
characteristic setting. While some may find the zero
characteristic setting more topologically natural or intuitive,
others may see this view as a restriction which, like many
restrictions, prevents the purest mathematical arguments from
being made. Indeed, after R.~Rentschler \cite{rentschler} provided
a description of algebraic $\mathbf{C}^+$-actions on
$\mathbf{C}^{2}$, M.~Miyanishi \cite{miyanishi1} demonstrated that
Rentschler's theorem extends naturally without any technical
difficulty to an affine plane of arbitrary characteristic, and
many other papers appeared concurrently which studied such
characteristic free techniques (e.g. the work of M.~Miyanishi and
Y.~Nakai \cite{miyanishi2,miynak,nakai}). Also, after T.~Fujita,
M.~Miyanishi and T.~Sugie \cite{fujita,miysug} affirmatively
solved the Zariski cancellation problem for the affine plane
$\mathbf{C}^2$, P. Russell \cite{russell} gave a simplified
treatment of their proof which erased the characteristic zero
restriction.

The purpose of the present paper is to place the AK invariant in a
characteristic free environment. We provide the definition and
basic ideas, and we demonstrate computational techniques, closely
following earlier treatments due to H.~Derksen, O.~Hadas, and
L.~Makar-Limanov \cite{derksen,DHML,ML1}, which gain no
simplification from the characteristic zero assumption. As an
illustration, we compute the AK invariant for $x + x^2 y + z^2 +
t^3 = 0$ over any field $\K$ and obtain the same result which was
originally found over $\mathbf{C}$. Similar efforts have been
successfully applied toward a generalization of the Zariski
cancellation problem \cite{CML}.

\section{Methods}
\subsection*{Exponential maps, the AK invariant, and locally finite iterative higher derivations}

Let $A$ be a $\K$-algebra. Suppose $\f:A \to A^{[1]}$ is a
$\K$-algebra homomorphism. We write $\f = \f_U:A \to A[U]$ if we
wish to emphasize an indeterminate $U$. We say that $\f$ is an
\emph{exponential map on $A$\/} if it satisfies the following two
additional properties.
\begin{itemize}
\item[(i)] $\varepsilon_0 \f_U$ is the identity on $A$, where
$\varepsilon_0:A[U] \to A$ is evaluation at $U=0$.\label{exp1}

\item[(ii)] $\f_S \f_U = \f_{S+U}$, where $\f_S$ is extended by
$\f_S(U)=U$ to a homomorphism $A[U] \to A[S,U]$.\label{exp2}
\end{itemize}
(When $A$ is the coordinate ring of an affine variety $\Spec(A)$
over $\K$, the exponential maps on $A$ correspond to algebraic
actions of the additive group $\K^+$ on $\Spec(A)$~\cite[\S
9.5]{essen}.) Define
\[
A^{\f} = \{ a \in A \,|\, \f(a)=a \},
\]
a subalgebra of $A$ called the \emph{ring of $\f$-invariants}. Let
$\EXP(A)$ denote the set of all exponential maps on $A$. We define
the \emph{AK invariant}, or \emph{ring of absolute constants of
$A$\/} as
\[
\ak(A) = \bigcap_{\f \in \EXP(A)} A^{\f}.
\]
This is a subalgebra of $A$ which is preserved by isomorphism.
Indeed, any isomorphism $f:A \to B$ of $\K$-algebras restricts to
an isomorphism $f: \ak(A) \to \ak(B)$. To understand this, observe
that if $\f \in \EXP(A)$ then $f \f f^{-1} \in \EXP(B)$. Remark
that $\ak(A) = A$ if and only if the only exponential map on $A$
is the standard inclusion $\f(a) = a$ for all $a \in A$.

\begin{example}
By considering exponential maps of the form $\f_i(X_j) = X_j +
\delta_{ij}U$, where $\delta_{ij}$ is the Kronecker delta, one can
see that $\ak(\K^{[n]}) = \K$ for each natural number $n$. When
$n=1$, this characterizes $\K^{[1]}$ (see Lemma~\ref{L:dim1}).
However, if $A$ is a domain with transcendence degree $n \geq 2$
over $\K$, then $\ak(A) = \K$ does not imply that $A \cong
\K^{[n]}$~\cite{bandman}.
\end{example}

It is often helpful to view a given $\f \in \EXP(A)$ as a sequence
in the following way. For each $a \in A$ and each natural number
$n$, let $D^n(a)$ denote the $U^n$-coefficient of $\f(a)$. Let $D
= \{ D^0, D^1, D^2, \ldots \}$. To say that $\f$ is a $\K$-algebra
homomorphism is equivalent to saying that the sequence
$\{D^i(a)\}$ has finitely many nonzero elements for each $a \in
A$, that $D^n:A \to A$ is $\K$-linear for each natural number $n$,
and that the Leibniz rule
\[
D^n(ab) = \sum_{i+j=n} D^i(a) D^j(b)
\]
holds for all natural numbers $n$ and all $a,b \in A$. The above
properties (i) and (ii) of the exponential map $\f$ translate into
the following properties of $D$.
\begin{itemize}
\item[(i')] $D^0$ is the identity map.

\item[(ii')] (iterative property) For all natural numbers $i,j$,
\[
D^i D^j = \binom{i+j}{i} D^{i+j}.
\]
\end{itemize}

Due to all of these properties, the collection $D$ is called a
\emph{locally finite iterative higher derivation on $A$}. More
generally, a \emph{higher derivation on $A$\/} is a collection $D
= \{D^i\}$ of $\K$-linear maps on $A$ such that $D^0$ is the
identity and the above Leibniz rule holds. The notion of higher
derivations is due to H.~Hasse and F.K.~Schmidt~\cite{hasse}. When
the characteristic of $A$ is zero, each $D^i$ is determined by
$D^1$, which is a locally nilpotent derivation on $A$. In this
case, $\f = \exp(U D^1) = \sum_i \frac{1}{i!}(U D^1)^i$ and
$A^{\f}$ is the kernel of $D^1$.

The above discussion of exponential maps, locally finite iterative
higher derivations, and the AK invariant makes sense more
generally for any (not necessarily commutative) ring. However, we
will not need this generality.

\subsection*{Degree functions and related lemmas}

Given an exponential map $\f: A \to A[U]$ on a domain $A$ over
$\K$, we can define the \emph{$\f$-degree\/} of an element $a \in
A$ by $\deg_{\f}(a) = \deg_U(\f(a))$ (where $\deg_U(0) = -
\infty$). Note that $A^{\f}$ consists of all elements of $A$ with
non-positive $\f$-degree. The function $\deg_{\f}$ is a degree
function on $A$, i.e. it satisfies these two properties for all
$a,b \in A$.
\begin{itemize}
\item[(i)] $\deg_{\f}(ab) = \deg_{\f}(a) + \deg_{\f}(b)$.

\item[(ii)] $\deg_{\f}(a+b) \leq \max \{\deg_{\f}(a),\deg_{\f}(b)
\}$.
\end{itemize}
Equipped with these notions, we now collect some useful facts.

\begin{lemma}\label{L:facts}
Let $\f$ be an exponential map on a domain $A$ over $\K$. Let $D =
\{D^i\}$ be the locally finite iterative higher derivation
associated to $\f$.

\renewcommand{\theenumi}{\alph{enumi}}
\begin{itemize}
\item[(a)] If $a,b \in A$ such that $ab \in A^{\f} \setminus 0$,
then $a,b \in A^{\f}$. In other words, $A^{\f}$ is factorially
closed in $A$.\label{factsa}

\item[(b)] $A^{\f}$ is algebraically closed in $A$.\label{factsb}

\item[(c)] For each $a \in A$, $\deg_{\f}(D^i(a)) \leq
\deg_{\f}(a) - i$. In particular, if $a \in A \setminus 0$ and $n
= \deg_{\f}(a)$, then $D^n(a) \in A^{\f}$.\label{factsc}
\end{itemize}
\end{lemma}

\begin{proof}
(a): We have $0 = \deg_{\f}(ab) = \deg_{\f}(a) + \deg_{\f}(b)$,
which implies that
$\deg_{\f}(a) = \deg_{\f}(b) = 0$.\\
(b): If $a \in A \setminus 0$ and $c_n a^n + \cdots + c_1 a + c_0
= 0$ is a polynomial relation with minimal possible degree $n \geq
1$, where each $c_i \in A^{\f}$ with $c_0 \ne 0$, then $a(c_n
a^{n-1} + \cdots +c_1) = -c_0 \in A^{\f}
\setminus 0$. By part (a), $a \in A^{\f}$.\\
(c): Use the iterative property of $D$ to check that
$D^j(D^i(a))=0$ whenever $j > \deg_{\f}(a) - i$.
\end{proof}

\begin{lemma}\label{L:facts2}
Let $\f$ be a nontrivial exponential map (i.e not the standard
inclusion) on a domain $A$ over $\K$ with $\ch(\K)=p \geq 0$. Let
$x \in A$ with minimal positive $\f$-degree $n$.

\renewcommand{\theenumi}{\alph{enumi}}
\begin{itemize}
\item[(a)] $D^i(x) \in A^{\f}$ for each $i \geq 1$. Moreover,
$D^i(x)=0$ whenever $i \geq 1$ is not a power of
$p$.\label{facts2a}

\item[(b)] If $a \in A \setminus 0$, then $n$ divides
$\deg_{\f}(a)$.\label{facts2b}

\item[(c)] Let $c = D^n(x)$. Then $A$ is a subalgebra of
$A^{\f}[c^{-1}][x]$, where $A^{\f}[c^{-1}] \subseteq
\Frac(A^{\f})$ is the localization of $A^{\f}$ at
$c$.\label{facts2c}

\item[(d)] Let $\trdeg_{\K}$ denote transcendence degree over
$\K$. If $\trdeg_{\K}(A)$ is finite, then $\trdeg_{\K}(A^{\f}) =
\trdeg_{\K}(A) - 1$.\label{facts2d}
\end{itemize}
\end{lemma}

\begin{proof}
In proving parts (a) and (b) we will utilize the following fact.
If $p$ is prime and $i = p^j q$ for some natural numbers $i,j,q$,
then $\binom{i}{p^j} \equiv q \pmod{p}$~\cite[Lemma
5.1]{isaacs}.\\
(a): By part (c) of Lemma~\ref{L:facts}, $D^i(x) \in A^{\f}$ for
all $i \geq 1$. If $p=0$ then $n=1$, for given any element in $A
\setminus A^{\f}$ we can find an element with $\f$-degree 1 by
applying the locally nilpotent derivation $D^1$ sufficiently many
times. In this case, the second statement is immediate. Suppose
now that $p$ is prime and that $i>1$ is not a power of $p$, say $i
= p^j q$, where $j$ is a nonnegative integer and $q \geq 2$ is an
integer not divisible by $p$. Then $D^{i - p^j}(x) \in A^{\f}$ and
\[
0 = D^{p^j} D^{i-p^j}(x) = \binom{i}{p^j}D^i(x) = q D^i(x).
\]
We can divide by $q$ to conclude that $D^i(x)=0$.\\
(b): Again if $p=0$ then $n=1$ and the claim is obvious. Assume
that $p$ is prime. By part (a) we have $n = p^m$ for some integer
$m \geq 0$. If $m=0$, the claim is immediate. Assume that $m>0$.
Let $d = \deg_{\f}(a)$. Suppose that $p$ does not divide $d$. By
part (c) of Lemma~\ref{L:facts}, $\deg_{\f}(D^{d-1}(a)) \leq 1$.
Now, $D^1D^{d-1}(a) = d D^d(a) \ne 0$. So $\deg_{\f}(D^{d-1}(a)) =
1 < n$, contradicting the minimality of $n$. Hence we can write $d
= p^k d_1$ with $k \geq 1$ and $d_1$ not divisible by $p$. Making
a similar computation, $D^{p^k} D^{d-p^k}(a) = d_1 D^d(a) \ne 0$.
This implies that $\deg_{\f}(D^{d-p^k}(a)) = p^k$. Since $n = p^m$
is minimal,
we must have $k \geq m$, and so $n$ divides $d$.\\
(c): Let $a \in A \setminus 0$. By part (b) we can write
$\deg_{\f}(a) = ln$ for some natural number $l$. If $l = 0$ then
$a \in A^{\f}$ and we are done. We use induction on $l>0$.
Elements $c^l a$ and $D^{ln}(a) x^l$ both have $\f$-degree $ln$.
Let us check that $D^{ln}(c^l a) = D^{ln}(D^{ln}(a) x^l)$. First,
$D^{ln}(c^l a) = c^l D^{ln}(a)$ by the Leibniz rule and because
$c^l$ is $\f$-invariant. Secondly, since $D^{ln}(x^l) = D^n(x)^l =
c^l$ and $D^{ln}(a)$ is $\f$-invariant, we see that
$D^{ln}(D^{ln}(a) x^l) = c^l D^{ln}(a)$ as well. (Remark: Though
the equality $D^{ln}(x^l) = D^{n}(x)^l$ does follow from the
Leibniz rule, it may be more immediately observed as follows.
$D^{n}(x)$ is the leading $U$-coefficient of $\f(x)$, and $\f$ is
a homomorphism. Hence the leading $U$-coefficient of $\f(x^l)$ is
also that of $\f(x)^l$.) Therefore, the element $y = c^l a -
D^{ln}(a)x^l$ has $\f$-degree less than $ln$ and hence less than
or equal to $(l-1)n$. By the inductive hypothesis, $y \in
A^{\f}[c^{-1}][x]$. So $a = c^{-l}(y
+ D^{ln}(a) x^l) \in A^{\f}[c^{-1}][x]$.\\
(d): This is immediate from part (c), together with part (b) of
Lemma~\ref{L:facts} which states that $A^{\f}$ is algebraically
closed in $A$.
\end{proof}

\begin{lemma}\label{L:dim1}
Let $A$ be a domain over $\K$ with $\trdeg_{\K}(A)=1$. Then
$\ak(A) = \K$ if and only if $A \cong \K^{[1]}$. Otherwise,
$\ak(A)=A$.
\end{lemma}

\begin{proof}
To see that $\ak(\K[X]) = \K$ (where $X$ is an indeterminate),
observe that $\psi(X) = X + U$ defines an exponential map on
$\K[X]$ with ring of $\psi$-invariants $\K$. Suppose $\ak(A) \ne
A$. Let $\f \in \EXP(A)$ be nontrivial. Part (b) of
Lemma~\ref{L:facts} implies that $A^{\f} = \K$. By part (c) of
Lemma~\ref{L:facts2}, $A \subseteq \K[x]$ for some $x \in A$ with
minimal positive $\f$-degree. So $A = \K[x]$.
\end{proof}

Thus $\K^{[1]}$ is the only trancendence degree 1 domain over $\K$
which admits nontrivial exponential maps.

\begin{lemma}\label{L:hopital}
Let $\f$ be an exponential map on a domain $A$ over $\K$. Extend
$\f$ to a homomorphism $\f: \Frac(A) \to \Frac(A)(U)$ by
$\f(ab^{-1}) = \f(a)\f(b)^{-1}$, and let $\Frac(A)^{\f} = \{f \in
\Frac(A) \, | \, \f(f)=f \}$. Then $\Frac(A)^{\f} =
\Frac(A^{\f})$.
\end{lemma}

\begin{proof}
It is clear that $\Frac(A^{\f}) \subseteq \Frac(A)^{\f}$. Suppose
that $a,b \in A \setminus 0$, such that $ab^{-1} \in
\Frac(A)^{\f}$. Then
\[
ab^{-1} = \f(ab^{-1}) = \f(a)\f(b)^{-1} = \left(\sum_i D^i(a)U^i
\right) \left( \sum_i D^i(b)U^i \right)^{-1},
\]
and so
\[
\sum_i aD^i(b)U^i = \sum_i bD^i(a)U^i.
\]
Compare the leading coefficients to see that $a$ and $b$ have the
same $\f$-degree, say $n$, and $a b^{-1} = D^n(a) D^n(b)^{-1} \in
\Frac(A^{\f})$.
\end{proof}

Let $A$ be a domain over $\K$, and let $\f \in \EXP(A)$. Let $B$
be the domain over $\Frac(A^{\f})$ obtained by localizing $A$ at
the multiplicative set $A^{\f} \setminus 0$. By
Lemma~\ref{L:hopital}, $\f$ extends to a
$\Frac(A^{\f})$-homomorphism $\f:B \to B[U]$ which is an
exponential map on $B$ with ring of invariants $\Frac(A^{\f})$.
\subsection*{Homogenization of an exponential map}

Let $A$ be a domain over $\K$. Let $\mathbf{Z}$ denote the
integers. Suppose that $A$ has a $\mathbf{Z}$-filtration $\{ A_n
\}$. This means that $A$ is the union of linear subspaces $A_n$
with these properties: $A_i \subseteq A_j$ whenever $i \leq j$,
and $A_i \cdot A_j \subseteq A_{i+j}$ for all $i,j \in
\mathbf{Z}$. Additionally, suppose that
\[
(A_i \setminus A_{i-1}) \cdot (A_j \setminus A_{j-1}) \subseteq
A_{i+j} \setminus A_{i+j-1}
\]
for all $i,j \in \mathbf{Z}$. This will be the case if the
filtration is induced by a degree function. Suppose also that
$\chi$ is a set of generators for $A$ over $\K$ with the following
property: if $a \in A_i \setminus A_{i-1}$ then we can write $a =
\sum_{I} c_{I} \mathbf{x}^{I}$, a summation of monomials $c_{I}
\mathbf{x}^{I}$ built from $\chi$ which are all contained in
$A_i$. This is not an unreasonable property. It merely asserts
some homogeneity on the generating set $\chi$.

Given $a \in A \setminus 0$ there exists $i \in \mathbf{Z}$ for
which $a \in A_i \setminus A_{i-1}$. Write
\[
\tp{a} = a + A_{i-1} \in A_i / A_{i-1},
\]
the \emph{top part of $a$}. We can construct a graded $\K$-algebra
\[
\gr(A) = \bigoplus_{n \in \mathbf{Z}} A_n / A_{n-1}.
\]
Addition on $\gr(A)$ is given by its vector space structure. Given
$\tp{a} = a + A_{i-1}$ and $\tp{b} = b + A_{j-1}$, define $\tp{a}
\, \tp{b} = ab + A_{i+j-1}$. Note that $\tp{a} \, \tp{b} =
\tp{ab}$. Extend this multiplication to all of $\gr(A)$ by the
distributive law. By our assumption on the filtration, $\gr(A)$ is
a domain. Also, $\gr(A)$ is generated by the top parts of the
elements of $\chi$. Therefore, if $\chi$ is a finite set then
$\gr(A)$ is an affine domain.

Let $\grdeg$ be the degree function induced by the grading on
$\gr(A)$. By assigning a value to $\grdeg(U)$ for an indeterminate
$U$, we can extend the grading on $\gr(A)$ to $\gr(A)[U]$. Given
an exponential map $\f: A \to A[U]$ on $A$, the goal is to obtain
an exponential map $\tp{\f}$ on $\gr(A)$. For $a \in A$, let
$\grdeg(a)$ denote $\grdeg(\tp{a})$. Note that $\grdeg(\tp{a}) =
i$ if and only if $a \in A_i \setminus A_{i-1}$. Consequently,
$\grdeg$ can also be viewed as a degree function on $A$ and on
$A[U]$ once the value of $\grdeg(U)$ is determined.

Define
\begin{equation}\label{E:grdegU}
\grdeg(U) = \min \left\{ \frac{\grdeg(x)-\grdeg(D^{i}(x))}{i}
\biggm | x \in \chi, i \in \mathbf{Z}^+ \right\}.\tag{$\star$}
\end{equation}
Let us assume now that $\grdeg(U)$ does exist, i.e. is a rational
number. This will indeed occur whenever $\chi$ is a finite set, as
will be the case with our example of interest. If $x \in \chi$ and
$n$ is a natural number, then $\grdeg(D^n(x)U^n) \leq \grdeg(x)$
by our choice of $\grdeg(U)$. From this it follows by
straightforward calculation that $\grdeg(D^n(a) U^n) \leq
\grdeg(a)$ for all $a \in A$ and all natural numbers $n$. (Here we
use the homogeneity assumption imposed on $\chi$.) The reader can
easily work out the details or refer to \cite{thesis}. Note that
this inequality is sharp since
\[
\grdeg(U) = \frac{1}{n}(\grdeg(x) - \grdeg(D^n(x)))
\]
for some $x \in \chi$ and some positive integer $n$ (and also
since $D^0(a) = a$ for all $a \in A$).

For $a \in A$, let
\[
S(a) = \{n \, | \, \grdeg(D^n(a)) + n \grdeg(U) = \grdeg(a) \}.
\]
Define
\[
\tp{\f} (\tp{a}) = \sum_{n \in S(a)} \tp{D^n(a)} U^n
\]
and extend this linearly to define $\tp{\f} : \gr(A) \to
\gr(A)[U]$, the \emph{homogenization\/} or \emph{top part of
$\f$}. One can verify that $\tp{\f}$ is an exponential map on
$\gr(A)$. Refer to~\cite{DHML} for the case $A = \K^{[n]}$. The
proof of the general case is symbolically identical. Let
$\tp{A^{\f}}$ denote the domain generated by the top parts of all
elements in $A^{\f}$. The end result is

\begin{theorem}[H.~Derksen, O.~Hadas, L.~Makar-Limanov~\cite{DHML}]\label{T:toppart}
Let $A$ be a domain over $\K$ with $\mathbf{Z}$-filtration $\{ A_n
\}$ such that $(A_i \setminus A_{i-1}) \cdot (A_j \setminus
A_{j-1}) \subseteq A_{i+j} \setminus A_{i+j-1}$ for all $i,j \in
\mathbf{Z}$. Let $\f$ be a nontrivial exponential map on $A$.
Assume that $\grdeg(U)$ exists as defined above. Then $\tp{\f}$ as
defined above is a nontrivial exponential map on $\gr(A)$.
Moreover, $\tp{A^{\f}}$ is contained in $\gr(A)^{\tp{\f}}$.
\end{theorem}

An important special case of homogenization is when $A$ itself is
graded. Then we can filter $A$ so that $\gr(A)$ is canonically
isomorphic to $A$, and we can choose $\chi$ to be a set of
homogeneous generators of $A$. In this case the top part of $\f$
is a nontrivial exponential map on $A$ (assuming $\grdeg(U)$
exists).

\begin{example}\label{E:toppart} Let $A = \K[X,Y]$, where $\ch(\K) = p$, prime.
Define $\f \in \EXP(A)$ by $\f(X) = X$ and $\f(Y) = Y + U + X
U^p$. We can grade $A$ by setting $\grdeg(X)=\alpha$ and
$\grdeg(Y)=\beta$ (with $\grdeg(\lambda)=0$ for all $\lambda \in
\K^*$, and $\grdeg(0)= -\infty$). Since $\grdeg(D^i(X))= - \infty$
for all $i \geq 1$, $X$ will not contribute to the value of
$\grdeg(U)$. Therefore,

\begin{align*}
\grdeg(U) &= \min \left\{ \frac{\grdeg(Y)-\grdeg(1)}{1},
\frac{\grdeg(Y)-\grdeg(X)}{p} \right\}\\
          &= \min \left\{ \beta, \frac{\beta - \alpha}{p}
          \right\}.
\end{align*}
In any case, $\tp{\f}(X) = X$. If $\beta < \frac{1}{p}(\beta -
\alpha)$ then $\grdeg(U) = \beta$ and $\tp{\f}(Y) = Y + U$. If
$\beta = \frac{1}{p}(\beta - \alpha)$ then $\grdeg(U) = \beta$ and
$\tp{\f}(Y) = \f(Y)$. If $\beta > \frac{1}{p}(\beta - \alpha)$
then $\grdeg(U) = \frac{1}{p}(\beta - \alpha)$ and $\tp{\f}(Y) = Y
+ X U^p$.
\end{example}

\section{The Russell hypersurface}

Let $R = \K[X,Y,Z,T] / (X + X^2 Y + Z^2 + T^3)$. If we wish to
emphasize a choice of $\K$, we write $R = R_{\K}$. Let $x,y,z,t
\in R$ denote the cosets of $X,Y,Z,T$, respectively. We shall
prove

\begin{theorem}\label{T:hypersurface}
Suppose $\K$ is algebraically closed.  Then $\ak(R) = \K[x]$.
\end{theorem}

If $\K$ is any field, then $\ak(\K^{[3]}) = \K$. Also, if
$\tp{\K}$ is an algebraic closure of $\K$, then $R_{\K}
\otimes_{\K} \tp{\K} = R_{\tp{\K}}$. These observations lead
immediately to

\begin{corollary}
$R \ncong \K^{[3]}$ for any $\K$.
\end{corollary}

To prove Theorem~\ref{T:hypersurface}, let us start with
\begin{lemma}\label{L:technical1}
Suppose $\K$ is algebraically closed.  Then $x \in \ak(R)$.
\end{lemma}

\begin{proof}
Suppose that $\f: R \to R[U]$ is a nontrivial exponential map on
$R$. We want to show that $x \in R^{\f}$. Let us consider $R$ as a
subalgebra of $\K[x,x^{-1},z,t]$ with
\[
y = -x^{-2}(x + z^2 +t^3).
\]
Introduce a degree function $w_1$ on $\K[x,x^{-1},z,t]$ by
\begin{align*}
w_1(x)&=-1,\\
w_1(z)&=0,\\
w_1(t)&=0
\end{align*}
(with $w_1(\lambda)=0$ for all $\lambda \in \K^*$ and $w_1(0)= -
\infty$). Then
\[
w_1(y)=2.
\]
This induces a $\mathbf{Z}$-filtration $\{R_i\}$ on $R$, where
$R_i$ consists of all $r \in R$ with $w_1(r) \leq i$. Passing to
top parts, $\tp{y} = - \tp{x}^{\,-2}(\tp{z}^{\,2} +
\tp{t}^{\,3})$. So the corresponding graded domain $\gr(R)$ is
generated by $\tp{x},\tp{y},\tp{z},\tp{t}$ and subject to the
relation $\tp{x}^{\,2} \tp{y} + \tp{z}^{\,2} + \tp{t}^{\,3} = 0$.
Let us write $x,y,z,t$ in place of $\tp{x},\tp{y},\tp{z},\tp{t}$,
respectively. Then
\[
\gr(R) = \K[x,y,z,t]/(x^2 y + z^2 + t^3).
\]

For a first step we show
\begin{sublemma}
$R^{\f} \subseteq \K[x,z,t]$.
\end{sublemma}

\begin{proof}
Suppose that $f \in R^{\f}$ and $f \notin \K[x,z,t]$. It is clear
that the value $\grdeg(U)$ as defined by formula (\ref{E:grdegU})
exists for our filtration, since $R$ is finitely generated by
$\chi = \{x,y,z,t\}$. By Theorem~\ref{T:toppart}, $\f$ induces a
nontrivial exponential map $\tp{\f}$ on $\gr(R)$ with $\tp{f} \in
\gr(R)^{\tp{\f}}$. We can write
\[
\tp{f} = x^a y^b g(z,t)
\]
for some natural numbers $a$ and $b$ and some polynomial $g(z,t)
\in \K[z,t]$. By our assumption on $f$, we know $b$ must be
positive. We can assume that $a = 0$ or $a=1$, since a factor $x^2
y$ of $\tp{f}$ can be absorbed into $g(z,t)$ by the relation on
$\gr(R)$. If $a=1$, then we can replace $f$ by $f^2$, and so we
may assume that $a=0$. So now $\tp{f} = y^b g(z,t)$. Also, since
$\trdeg_{\K}(R^{\f}) = 2$ (part (d) of Lemma~\ref{L:facts2}), we
can assume that $g(z,t)$ is not a constant polynomial by replacing
$f$ if necessary.

Since $\gr(R)^{\tp{\f}}$ is factorially closed (part (a) of
Lemma~\ref{L:facts}), both $y$ and $g(z,t)$ belong to
$\gr(R)^{\tp{\f}}$. Let us introduce a new grading on $\gr(R)$ by
$w_2(x)=6$, $w_2(y)=-6$, $w_2(z)=3$, and $w_2(t)=2$. The
corresponding graded domain (still call it $\gr(R)$) is again
isomorphic to $\K[x,y,z,t]/(x^2 y + z^2 + t^3)$, and let us
continue to write $x,y,z,t$ in place of
$\tp{x},\tp{y},\tp{z},\tp{t}$. Under this new grading, we can
write
\[
g(z,t) = \lambda z^n t^m \prod_{i} (z^2 + \mu_{i} t^3)
\]
for some natural numbers $n$ and $m$ and for some $\lambda, \mu_i
\in \K^*$.

The next step is to show that neither $z$ nor $t$ can be
$\tp{\f}$-invariant. Suppose that $z \in \gr(R)^{\tp{\f}}$. Then
$\gr(R)^{\tp{\f}} = \K[y,z]$. By virtue of the remark following
Lemma~\ref{L:hopital}, $\tp{\f}$ induces a nontrivial exponential
map on the domain
\[
\K(y,z)[x,t]/(x^2 y + z^2 + t^3).
\]
But this domain has transcendence degree 1 over $\K(y,z)$ and is
not isomorphic to $\K(y,z)^{[1]}$. This contradicts
Lemma~\ref{L:dim1}. Thus $z \notin \gr(R)^{\tp{\f}}$. Suppose now
that $t \in \gr(R)^{\tp{\f}}$. Then $\gr(R)^{\tp{\f}} = \K[y,t]$
and just as with $z$ we can obtain a nontrivial exponential map on
the domain
\[
\K(y,t)[x,z]/(x^2 y + z^2 + t^3).
\]
Again, this domain has transcendence degree 1 over $\K(y,t)$ but
is not isomorphic to $\K(y,t)^{[1]}$, contradicting
Lemma~\ref{L:dim1}. So $t \notin \gr(R)^{\tp{\f}}$ as well. Since
$\gr(R)^{\tp{\f}}$ is factorially closed, $n=0$ and $m=0$ in the
above factorization of $g(z,t)$.

It remains to consider the factors of $g(z,t)$ of the form $z^2 +
\mu_{i} t^3$. Ignoring multiplicity, there can only be one such
factor. For given two factors of this type, their difference
belongs to $\gr(R)^{\tp{\f}}$ from which we conclude that both
$z^2$ and $t^3$ belong to $\gr(R)^{\tp{\f}}$, which we have just
shown to be impossible. Also, $z^2 + t^3$ cannot be a factor of
$g(z,t)$, since otherwise we would have $x^2 y = -(z^2 + t^3) \in
\gr(R)^{\tp{\f}}$, which in turn implies that $\gr(R)^{\tp{\f}}=
\gr(R)$, contradicting the nontriviality of $\tp{\f}$. We can
therefore write $g(z,t)$ as
\[
g(z,t) = \lambda (z^2 + \mu t^3)^k
\]
for some positive integer $k$ and some $\lambda, \mu \in \K^*$
with $\mu \ne 1$. We will use the same trick that worked for $z$
and $t$. Let $S$ be the domain which results from localizing
$\gr(R)$ at $\gr(R)^{\tp{\f}} \setminus 0$. Recall that $y$ and
$z^2 + \mu t^3$ belong to $\gr(R)^{\tp{\f}}$. Note that we can
rewrite $x^2 y + z^2 + t^3=0$ as the relation
\[
x^2 y + (1 - \mu)t^3 + (z^2 + \mu t^3) = 0
\]
over $\Frac(\gr(R)^{\tp{\f}})$. From this we can see that $S$ has
transcendence degree 1 over $\Frac(\gr(R)^{\tp{\f}})$ but is not
isomorphic to $\Frac(\gr(R)^{\tp{\f}})^{[1]}$. We can extend
$\tp{\f}$ to a nontrivial exponential map on $S$ over
$\Frac(\gr(R)^{\tp{\f}})$ in the manner described after
Lemma~\ref{L:hopital}. This once again contradicts
Lemma~\ref{L:dim1}. We have now exhausted all possibilities. To
avoid a contradiction, we must have $f \in \K[x,z,t]$. So $R^{\f}
\subseteq \K[x,z,t]$.
\end{proof}

Continuing with the proof of Lemma~\ref{L:technical1}, we are now
in position to show that $x \in R^{\f}$. Suppose that it is not
the case. If $f \in R^{\f} \subseteq \K[x,z,t]$, write $f = x
f_1(x,z,t) + f_2(z,t)$. Then $f_2 \ne 0$ since $x \notin R^{\f}$.
Again consider $\gr(R)$ given by $w_1$. Now $w_1(x f_1(x, z,t))$
is negative, while $w_1(f_2(z,t))=0$, and so $\tp{f} = f_2(z,t)
\in \gr(R)^{\tp{\f}}$. Let $g \in R^{\f}$ be algebraically
independent with $f$ over $\K$. (Recall that
$\trdeg_{\K}(R^{\f})=2$ by part (d) of Lemma~\ref{L:facts2}.) We
write $g = xg_1(x,z,t) + g_2(z,t)$, where $0 \ne g_2(z,t) =
\tp{g}$. Suppose that $f_2$ and $g_2$ are algebraically dependent
over $\K$, say $P(f_2,g_2) = 0$. Then $P(f,g)$ is a nonzero
element of $R^{\f}$, but $P(f,g)$ is divisible by $x$. This
implies that $x \in R^{\f}$, contrary to our assumption. Hence it
must be that $f_2$ and $g_2$ are algebraically independent over
$\K$. Thus $\gr(R)^{\tp{\f}}$ contains two algebraically
independent elements of $\K[z,t]$. Since $\gr(R)^{\tp{\f}}$ is
algebraically closed in $\gr(R)$ (part (b) of
Lemma~\ref{L:facts}), we deduce that $\gr(R)^{\tp{\f}} = \K[z,t]$.
Now $x^2 y = -z^2 - t^3 \in \gr(R)^{\tp{\f}}$, and this implies
that $x,y \in \gr(R)^{\tp{\f}}$. But then $\tp{\f}$ is trivial.
This contradicts our assumption that $x \notin R^{\f}$. So $x \in
R^{\f}$ for every $\f \in \EXP(R)$, and the lemma is finally
proved.
\end{proof}

\begin{proof}[Proof of Theorem~\ref{T:hypersurface}]
We know that $\K[x] \subseteq \ak(R)$. To show the reverse
containment, we consider two maps $\f_1$ and $\f_2$ which one
easily verifies to be exponential maps on $R$. Define $\f_1 : R
\to R[U]$ by
\begin{align*}
\f_1(x) &=x,\\
\f_1(y) &=y + 2zU - x^2 U^2,\\
\f_1(z) &=z - x^2 U,\\
\f_1(t) &=t.
\end{align*}
The ring of $\f_1$-invariants is $\K[x,t]$. Define $\f_2 : R \to
R[U]$ by
\begin{align*}
\f_2(x) &=x,\\
\f_2(y) &=y + 3t^2 U - 3x^2 tU^2 + x^4 U^3,\\
\f_2(z) &=z,\\
\f_2(t) &=t - x^2 U.
\end{align*}
The ring of $\f_2$-invariants is $\K[x,z]$. So $\ak(R)$ is
contained in the intersection of these two rings, that being
$\K[x]$.
\end{proof}

As a final remark, L.~Makar-Limanov \cite{ML3} has recently taken
a simplified approach (again very similar to the proof given here)
to showing that $\ak(R_{\mathbf{C}}) = \mathbf{C}[x]$. This proof
uses the following fact.

\begin{lemma}[see \cite{ML3}]\label{L:0}
Let $A$ be a domain with characteristic zero. Let $n$ and $m$ be
natural numbers both at least 2. Let $\f \in \EXP(A)$ and $c_1,c_2
\in A^{\f} \setminus 0$. If $a,b \in A$ such that $c_1 a^n + c_2
b^m \in A^{\f} \setminus 0$, then $a,b \in A^{\f}$.
\end{lemma}

In fact, this lemma is still true when $A$ has prime
characteristic $p$, under the additional necessary hypothesis that
neither $n$ nor $m$ are powers of $p$ \cite{thesis}. This lemma
can then replace the many times that we used Lemma~\ref{L:hopital}
to contradict Lemma~\ref{L:dim1}. But because of the extra
assumption needed on Lemma~\ref{L:0}, this method fails when the
characteristic of $R$ is 2 or 3.

\end{document}